\documentclass[12pt,reqno]{amsart}

\usepackage{amsfonts,amsmath,amssymb}
\usepackage{graphicx}

\numberwithin{equation}{section}

\advance\textheight by \topskip
\newcommand{\gauss}[2]{\genfrac{[}{]}{0pt}{}{#1}{#2}}

\numberwithin{equation}{section}

\title[${q}$--identities of Fu and Lascoux]
{$\boldsymbol{q}$--identities of Fu and Lascoux proved by the
$\boldsymbol{q}$--Rice formula}

\author[H.~Prodinger]{Helmut Prodinger}
\address{Helmut Prodinger\\
The John Knopfmacher Centre for Applicable Analysis and Number Theory\\
 School of Mathematics\\
University of the Witwatersrand, Private Bag~3,
Wits, 2050 Johannesburg, South Africa}
\email{helmut@maths.wits.ac.za}

\date{April 8, 2004}

\begin{document}
\begin{abstract}
Two recent $q$--identities of Fu and Lascoux are proved by the $q$--Rice formula.
\end{abstract}

\maketitle

\section{Introduction}

Fu and Lascoux \cite{FuLa04}  (answering questions
of Corteel and Lovejoy \cite{CoLo04}) proved the following two identities:
\begin{multline}\label{Identity_1}
\sum_{i=1}^n\gauss ni(-1)^{i-1}(x+1)\dots(x+q^{i-1})\frac{q^{mi}}{(1-q^i)^m}\\
=\sum_{i=1}^n\big(1-(-x)^i\big)\frac{q^i}{1-q^i}\sum_{i\le i_2\le\dots\le i_m\le
n}
\frac{q^{i_2}}{1-q^{i_2}}\dots\frac{q^{i_m}}{1-q^{i_m}}
\end{multline}
and
\begin{multline}\label{Identity_2}
\sum_{i=0}^n\gauss ni(-1)^{i-1}(x+1)\dots(x+q^{i-1})\frac{q^i}{1-tq^i}
=-\frac{(q;q)_n}{(t;q)_{n+1}}\sum_{i=0}^n\frac{(t;q)_{i}}{(q;q)_i}(-xq)^i.
\end{multline}
Here, we use the usual notation $(x;q)_n=(1-x)(1-xq)\dots(1-xq^{n-1})$ and\newline
$\gauss nk=(q;q)_n/(q;q)_k(q;q)_{n-k}$, see \cite{Andrews76}.

In this short note, we will provide attractive proofs of these, using the
$q$--Rice formula, see \cite{Prodinger01amy} for some background and applications.

\section{Proof of Identity~\eqref{Identity_1}}

The $q$--Rice formula \cite{Prodinger01amy} allows to write an alternating sum as
a countour integral:
\begin{equation*}
\sum_{i=1}^n\gauss ni(-1)^{i-1}q^{\binom i2}f(q^{-i})=\frac1{2\pi i}
\int_{\mathcal C}\frac{(q;q)_n}{(z;q)_{n+1}}f(z)dz,
\end{equation*}
where the curve $\mathcal C$ encircles the poles $q^{-1},\dots,q^{-n}$ and no others.
For more technical details, see \cite{Prodinger01amy}. Under mild conditions, the
integral (and thus the sum) can be expressed as the negative sum of the further residues.
Thus, the computation of the alternating sum boils down to a residue computation. 

In our application, we must find $f(z)$ such that
\begin{align*}
f(q^{-i})&=(x+1)\dots(x+q^{i-1})\frac{q^{mi}}{(1-q^i)^m }q^{-\binom i2}\\
&=(1+x)\dots(1+\frac x{q^{i-1}})\frac{1}{(q^{-i}-1)^m }.
\end{align*}
Now
\begin{equation*}
(1+x)\dots(1+\frac x{q^{i-1}})
=\prod_{h\ge1}\frac{1+\frac{xq^h}{q^i}}{1+xq^h}.
\end{equation*}
and thus we take
\begin{equation*}
f(z)=\frac{1}{(z-1)^m }\prod_{h\ge1}\frac{1+xzq^h}{1+xq^h}.
\end{equation*}
The only extra pole is at $z=1$, and so the sum is given by
\begin{align*}
\textsf{SUM}&=-\text{Res}_{z=1}\frac{(q;q)_n}{(z;q)_{n+1}}
\frac{1}{(z-1)^m }\prod_{h\ge1}\frac{1+xzq^h}{1+xq^h}\\
&=-[(z-1)^{-1}]\frac{(q;q)_n}{(z;q)_{n+1}}
\frac{1}{(z-1)^m }\prod_{h\ge1}\frac{1+xzq^h}{1+xq^h}\\
&=[(z-1)^{m}]\frac{(q;q)_n}{(zq;q)_{n}}\prod_{h\ge1}\frac{1+xzq^h}{1+xq^h}\\
&=[w^m]\frac{1}{\big(1-w\frac{q}{1-q}\big)\dots\big(1-w\frac{q^n}{1-q^n}\big)}
\prod_{h\ge1} \bigg(1+\frac{xwq^h}{1+xq^h}\bigg).
\end{align*}

It is not hard to see that
\begin{align*}
\prod_{h\ge1} \bigg(1+\frac{xwq^h}{1+xq^h}\bigg)
=1-w\sum_{i\ge1}(-x)^i\frac{q^i}{1-q^i}\prod_{1\le h<i}\Big(1-w\frac{q^h}{1-q^h}
\Big).
\end{align*}
To sketch a proof, let us look at the coefficient of $w^2$:
\begin{align*}
\sum_{1\le h_1<h_2}&\frac{xq^{h_1}}{1+xq^{h_1}}\frac{xq^{h_2}}{1+xq^{h_2}}
=\sum_{1\le h_1<h_2,\;k_1\ge1}\frac{xq^{h_1}}{1+xq^{h_1}}
(-1)^{k_1-1}x^{k_1}q^{h_2k_1}\\
&=\sum_{1\le h_1,\;k_1\ge1}\frac{q^{h_1(k_1+1)}}{1+xq^{h_1}}
(-1)^{k_1-1}x^{k_1+1}\frac{q^{k_1}}{1-q^{k_1}}\\
&=\sum_{1\le h_1,\;k_1\ge1,\;k_2\ge0}{q^{h_1(k_1+1)}}q^{h_1k_2}
(-1)^{k_1+k_2-1}x^{k_1+k_2+1}\frac{q^{k_1}}{1-q^{k_1}}\\
&=\sum_{1\le h_1,\;1\le k_1<k_2}q^{h_1k_2}
(-1)^{k_2}x^{k_2}\frac{q^{k_1}}{1-q^{k_1}}\\
&=\sum_{1\le k_1<k_2}(-x)^{k_2}\frac{q^{k_1}}{1-q^{k_1}}\frac{q^{k_2}}{1-q^{k_2}}.
\end{align*}
If one does this, say, also for the coefficient of $w^3$, then one quickly discovers
the general pattern, and these coefficients are the same as the coefficients
of the right side. (Is there an easier proof?)

Now
\begin{equation*}
[w^m]\frac{1}{\big(1-w\frac{q}{1-q}\big)\dots\big(1-w\frac{q^n}{1-q^n}\big)}
=\sum_{i=1}^n\frac{q^i}{1-q^i}\sum_{i\le i_2\le\dots\le i_m\le
n}\frac{q^{i_2}}{1-q^{i_2}}\dots\frac{q^{i_m}}{1-q^{i_m}}
\end{equation*}
is already known (Dilcher's sum \cite{Dilcher95, Prodinger01amy}), so we are left to prove that
\begin{multline*}
\sum_{i=1}^{n\,(\infty)}(-x)^i\frac{q^i}{1-q^i}\sum_{i\le i_2\le\dots\le i_m\le
n}
\frac{q^{i_2}}{1-q^{i_2}}\dots\frac{q^{i_m}}{1-q^{i_m}}\\
=[w^m]\frac{1}{\big(1-w\frac{q}{1-q}\big)\dots\big(1-w\frac{q^n}{1-q^n}\big)}
w\sum_{i\ge1}(-x)^i\frac{q^i}{1-q^i}\prod_{1\le h<i}\Big(1-w\frac{q^h}{1-q^h}\Big)
\end{multline*}
In terms of generating functions, we should show that
\begin{multline*}
\sum_{i=1}^\infty(-x)^iwa_i\frac{1}{(1-wa_i)\dots(1-wa_n)}\\=
\frac{1}{(1-wa_1)\dots(1-wa_n)}
\sum_{i\ge1}(-x)^iwa_i\prod_{1\le h<i}(1-wa_h),
\end{multline*}
where we wrote $a_i=q^i/(1-q^i)$ (but it holds in general).
But this in equivalent to
\begin{equation*}
\sum_{i=1}^\infty(-x)^iwa_i\frac{(1-wa_1)\dots(1-wa_n)}{(1-wa_i)\dots(1-wa_n)}=
\sum_{i\ge1}(-x)^iwa_i\prod_{1\le h<i}(1-wa_h),
\end{equation*}
and thus proved.

\section{Proof of Identity~\eqref{Identity_2}}

This time we take
\begin{equation*}
f(z)=\frac{1}{z-t }\prod_{h\ge1}\frac{1+xzq^h}{1+xq^h}
\end{equation*}
and write
\begin{align*}
\textsf{SUM}=\frac1{2\pi i}\int_{\mathcal C}\frac{(q;q)_n}{(z;q)_{n+1}}\frac
{1}{z-t }\prod_{h\ge1}\frac{1+xzq^h}{1+xq^h}dz.
\end{align*}
Now we use Cauchy's formula:
\begin{equation*}
\prod_{h\ge1}\frac{1+xzq^h}{1+xq^h}=\frac{(-xzq;q)_\infty}{(-xq;q)_\infty}=
\sum_{m\ge0}\frac{(z;q)_m}{(q;q)_m}(-xq)^m.
\end{equation*}
However, for the residues at $z=q^{-i}$, $i=0,\dots,n$, only the terms for $m\le n$
are relevant. Henceforth we may write
\begin{align*}
\textsf{SUM}=\frac1{2\pi i}\int_{\mathcal C}\frac{(q;q)_n}{(z;q)_{n+1}}\frac
{1}{z-t }\sum_{m=0}^n\frac{(z;q)_m}{(q;q)_m}(-xq)^m dz.
\end{align*}
For outside residues, there is only one, at $z=t$, and therefore
\begin{align*}
\textsf{SUM}&=-\text{Res}_{z=t}\frac{(q;q)_n}{(z;q)_{n+1}}\frac
{1}{z-t }\sum_{m=0}^n\frac{(z;q)_m}{(q;q)_m}(-xq)^m\\
&=-\frac{(q;q)_n}{(t;q)_{n+1}}\sum_{m=0}^n\frac{(t;q)_m}{(q;q)_m}(-xq)^m.
\end{align*}
This is clearly equivalent to the formula of Fu and Lascoux.

\bibliographystyle{plain}


\end{document}